\input amstex
\documentstyle{amsppt}
\magnification=\magstep1
 \hsize 13cm \vsize 18.35cm \pageno=1
\loadbold \loadmsam
    \loadmsbm
    \UseAMSsymbols
\topmatter
\NoRunningHeads
\title
$\text{ A new Kim's type  Bernoulli and Euler Numbers and}$
$\text{ related identities and zeta and $L$-functions}$
\endtitle
\author
 Y. Simsek${}^{\ast}$, T. Kim and D. Kim
\endauthor
\abstract In  this paper, by using $q$-deformed bosonic $p$-adic
integral, we give  $\lambda$-Bernoulli numbers and polynomials,
 we prove Witt's type formula of  $\lambda$-Bernoulli polynomials
 and  Gauss multiplicative formula for $\lambda$-Bernoulli
 polynomials.
 By using derivative operator to the generating functions of $\lambda$-Bernoulli polynomials
 and generalized $\lambda$-Bernoulli numbers,
we give Hurwitz type $\lambda$-zeta functions and Dirichlet's type
  $\lambda$-$L$-functions; which are interpolated $\lambda$-Bernoulli
 polynomials and generalized $\lambda$-Bernoulli numbers,
 respectively. We give generating function  of $\lambda$-Bernoulli
 numbers with order $r$. By using Mellin transforms to their function,
   we prove relations between multiply zeta
 function and $\lambda$-Bernoulli polynomials  and ordinary
 Bernoulli numbers of order $r$ and $\lambda$-Bernoulli numbers,
 respectively. We also study on  $\lambda$-Bernoulli numbers and
 polynomials in the space of locally constant. Moreover, we define
 $\lambda$-partial   zeta function  and interpolation function.
\endabstract
\thanks  2000 AMS Subject Classification: 11S80, 11B68, 11M99,
32D30
\newline $\ast$ This paper was  supported by the Scientific Research Project Adminstration  Akdeniz University.
\endthanks
\endtopmatter

\document

{\bf\centerline {\S 0. Introduction, definitions and notations}}

 \vskip 20pt

Throughout this paper, $\Bbb Z$, $\Bbb Z_p$,  $\Bbb Q_p$ and $\Bbb
C_p$ will be denoted by the ring of rational integers, the ring of
$p$-adic  integers,  the field of $p-$adic rational numbers and
the completion of the algebraic closure of $\Bbb Q_p$,
respectively.
 Let $\nu_p$ be
the normalized exponential  valuation of $\Bbb C_p$ with $|p|_p
=p^{-\nu_p (p)} =\frac{1}{p}$, cf. [2,3,4,5,6,7,8,9,16,17, 21,27].

When one talks of $q-$extension, $q$ considered in many ways such
as an indeterminate, a complex number $q\in\Bbb C$, as $p$-adic
number $q\in\Bbb C_p$. If $q\in\Bbb C$ one normally assumes that
$\vert q \vert <1$. If $q\in\Bbb C_p$, we normally assume that
$\vert q-1 \vert_p < p^{-\frac{1}{p-1}}$ so that $q^x =\exp(x \log
q )$ for $\vert x \vert_p \leq 1$. We use the following notations:
$$
[x] =[x:q] = \dfrac{1-q^x }{ 1-q}, \quad \text{ cf
}[3,4,5,6,8,9,24,25,27].
$$
Observe that  when $\lim_{q\to 1} [x] =x$, for any  $x$ with
$\vert
 x \vert_p \leq 1$ in the present $p$-adic case $[x: a] =\dfrac{1-a^x }{
 1-a}$.

 Let $d$ be a fixed integer and let  $p$ be a fixed prime number.
 For any positive integer $N$, we set
 $$\split
& \Bbb X = \lim_{\overleftarrow{N} } \left( \Bbb Z/ dp^N \Bbb Z
\right) ,\cr & \Bbb X^\ast = \underset {{0<a<d p}\atop
{(a,p)=1}}\to {\cup} (a+ dp \Bbb Z_p ), \cr & a+d p^N \Bbb Z_p =\{
x\in \Bbb X | x \equiv a \pmod{dp^n}\},\endsplit$$ where $a\in
\Bbb Z$ lies in $0\leq a < d p^N$. We assume that $u\in\Bbb C_p$
with $\vert 1-u \vert_p \geq 1$. cf. [3,4,5,6,7,8,24, 27].

For $x\in\Bbb Z_p$, we say that $g$ is a uniformly differentiable
function at point $a\in \Bbb Z_p$, and write $g\in UD(\Bbb Z_P )$,
the set of uniformly differentiable functions, if the difference
quotients,
$$
F_g (x,y )= \dfrac{g(y)-g(x)}{ y-x},
$$
have a limit $l =g^\prime (a)$ as $(x,y) \to (a,a)$. For $f\in
UD(\Bbb Z_p)$,  the {\it $q$-deformed bosonic $p$-adic integral}
was defined as
$$\eqalignno{
I_q (f) &=\int_{\Bbb Z_p} f(x) d \mu_q (x) \cr &=\lim_{N \to
\infty } \sum_{x=0}^{ p^N -1} f(x)\mu_q(x +p^N \Bbb Z_p )&(A) \cr
& =\lim_{N \to \infty } \sum_{x=0}^{ p^N -1} f(x)
\dfrac{q^x}{[p^N]} , \text{ [4,5,9].}}
$$
By Eq-(A), we have
$$
\lim_{q \to -q}  I_q (f)= I_{-q} (f) = \int_{\Bbb Z_p} f(x) d
\mu_{-q } (x).
$$
This integral, $I_{-q} (f)$, give the $q$-deformed integral
expression of fermioinc.  The classical Euler numbers were defined
by means of the following  generating function:
$$
\dfrac{2}{ e^t +1} = \sum_{m=0}^\infty E_m \dfrac{ t^m }{ m!} ,
\quad \vert t \vert < 2\pi . \quad [6,7,19,21].
$$

Let $u$ be algebraic in complex number field. Then Frobenius-Euler
polynomials [6,7,19,21] were defined by
$$
\eqalignno{ \dfrac{1-u}{e^t -u} e^{xt} &= e^{H(u,x) t}=
\sum_{m=0}^{\infty} H_m (u,x )\dfrac{ t^m }{ m!}  ,&(A1)}
$$
where we use technical method's notation by replacing  $H^m (u,
x)$ by $H_m (u, x)$ symbolically.  In case $x=0$, $H_m (u, 0) =H_m
(u)$, which is called Frobenius-Euler number. The Frobenius-Euler
polynomials of order $r$, denoted by $H_n^{(r)}(u,x )$, were
defined by
$$
\left( \dfrac{1-u}{ e^t -u} \right)^r  e^{tx} = \sum_{n=0}^\infty
H_n^{(r)} (u,x) \dfrac{t^n}{n!} \quad \text{ cf.} [7,10,26,27].
$$
The values at $x=0$ are called Frobenius-Euler numbers of order
$r$. When $r=1$, these numbers and polynomials are reduced to
ordinary Frobenius-Euler numbers and polynomials. In the usual
notation, the $n$-th Bernoulli polynomial were defined by means of
the following generating function:
$$
\left( \dfrac{t}{ e^t -1} \right)  e^{tx} =  \sum_{n=0}^\infty B_n
(x)  \dfrac{t^n}{n!}.
$$

For $x=0$, $B_n (0) =B_n$ are said to be the $n$-th Bernoulli
numbers. The Bernoulli polynomials of order $r$ were defined by

$$
\left( \dfrac{t}{ e^t -1} \right)^r  e^{tx} =  \sum_{n=0}^\infty
B_n^{(r)} (x)  \dfrac{t^n}{n!}
$$
and $B_n^{(r)} (0) =B_n^{(r)}$ are called  the  Bernoulli numbers
or order $r$. Let $x, w_1 , w_2 , \cdots,w_r $ be complex numbers
with positive real parts. When the generalized Bernoulli numbers
and polynomials were defined by means of the following generating
function:
$$
\dfrac{w_1 w_2 \cdots w_r t^r e^{xt} }{ (e^{w_1 t} -1)(e^{w_2 t}
-1)\cdots (e^{w_r t} -1)} = \sum_{n=0}^\infty B_n^{(r)} (x ~| ~
w_1, w_2 \cdots, w_r ) \dfrac{t^n}{n!}$$ and $B_n^{(r)} (0 ~| ~
w_1, w_2 \cdots, w_r )= B_n^{(r)} ( w_1, w_2 \cdots, w_r )$, cf.
[13,15].

The Hurwitz zeta function is defined by $$\zeta (s, x)
=\sum_{n=0}^\infty \dfrac{1}{(x+n)^s},$$ $\zeta(s,1)= \zeta(s)$,
which is Riemann zeta function. The multiple zeta functions
[12,27] were defined by
$$\eqalignno{ &
\zeta_r (s) = \sum_{0<n_1 < n_2 <\cdots< n_r} \dfrac{1}{(n_1
+\cdots + n_r )^s }.&(C)}
$$

We summarize our paper as follows:

In section 1, by using $q$-deformed bosonic $p$-adic integral,
generating function of $\lambda$-Bernoulli numbers and polynomials
are given. We obtain many new identities related to these numbers
and polynomials. We proved Gauss multiplicative formula for
$\lambda$-Bernoulli numbers.  Witt's type formula of
$\lambda$-Bernoulli polynomials is given.

 In section 2, by using
$\left. \left( \dfrac{d}{dt} \right)^k \right\vert_{t=0}$
derivative operator to the generating function of the
$\lambda$-Bernoulli numbers, we define new relations and Hurwitz'
type $\lambda$-function, which interpolates $\lambda$-Bernoulli
polynomials at negative integers.

In section 3, by using same method of  section 2, we give
Dirichlet type $\lambda$-$L$-function which interpolates
generalized $\lambda$-Bernoulli numbers.

 In section 4, generating
function of $\lambda$-Bernoulli numbers of order $r$ is defined,
by using Cauchy residue theorem and Mellin transforms to this
function, we proved relation between multiple zeta function and
$\lambda$-Bernoulli numbers of order $r$.

In section 5, we give some important identities related to
generalized $\lambda$-Bernoulli numbers of order $r$.

In section 6, we study on $\lambda$-Bernoulli numbers and
polynomials in the space of locally constant. In this section, we
also define $\lambda$-partial zeta function which interpolates
$\lambda$-Bernoulli numbers at negative integers.

In section 7, we give $p$-adic interpolation functions.

 \vskip 20pt

{\bf\centerline {\S 1. $\lambda$-Bernoulli numbers}}

 \vskip 20pt

In this section, by using Eq-(A), we give integral equation of
bosonic $p$-adic integral. By using this integral equation we
define generating function of  $\lambda$-Bernoulli polynomials. We
give fundamental properties of the $\lambda$-Bernoulli numbers and
polynomials. We also give some new identities related to
$\lambda$-Bernoulli numbers and polynomials.  We prove Gauss
multiplicative formula for $\lambda$-Bernoulli numbers as well.
Witt's type formula of $\lambda$-Bernoulli polynomials is given.

To give the expression of bosonic $p$-adic integral in Eq-(A), we
consider the limit
$$\eqalignno{ &
I_1 (f) =\lim_{q\to 1} I_q (f)= \int_{\Bbb Z_p} f(x) d \mu_1
(x).&(0)}
$$
Bosonic $p$-adic integral on $\Bbb Z_p$(= $p$-adic invariant
integral on $\Bbb Z_p$)
$$\eqalignno{ &
I_1 (f_1 ) = I_1 (f) + f^\prime (0),&(1)}
$$
where $f_1 (x) =f(x+1)$, integral equation for bosonic $p$-adic
integral. Let $ C_{p^n }$ be the space of primitive $p^n$-th root
of unity, $$ C_{p^n}= \{ \zeta ~ | ~ \zeta^{p^n} =1\}.$$

Then, we denote
$$
T_p =\lim_{n \to \infty} C_{p^n} =\underset{n\geq 0} \to \cup
C_{p^n}.
$$

For $\lambda \in\Bbb Z_p$, we take $f(x) = \lambda^x e^{tx}$, and
$f_1 (x) = e^{t} \lambda f(x)$. Thus we have
$$ \eqalignno{ &
f_1 (x) -f(x) =(\lambda e^t -1)f(x). &(2)} $$

By substituting (2) into (1), we get
$$ \eqalignno{ &
(\lambda e^t -1)I_1 (f) = f^\prime (0). &(2a)}
$$
Consequently, we have

$$ \eqalignno{ &  \dfrac{\log \lambda +t}{\lambda e^t -1} :=
\sum_{n=0}^\infty B_n (\lambda) \dfrac{t^n}{ n!}  . &(3)}
$$
By using Eq-(3), we obtain
$$
\lambda (B(\lambda ) +1 )^n - B_n (\lambda ) =  \cases
\log \lambda , & \text{if\   $n=0$}\\
1, & \text{if \  $n=1$}\\
0, & \text{if \  $n>1$,}\\
\endcases
$$
with the usual convention of replacing $B_n (\lambda )$ by $B^n
(\lambda)$. We give some $B_n (\lambda)$ numbers as follows:
$$
B_0 (\lambda) =\dfrac{ \log \lambda}{\lambda -1}, \quad B_1
(\lambda) =\dfrac{ \lambda -1-\lambda \log \lambda}{(\lambda
-1)^2}, \cdots .
$$

We note that, if $\lambda \in T_p$, for some $n\in\Bbb N$, then
Eq-(2a) is reduced to the following generating function:
$$\eqalignno{ &
\dfrac{t}{\lambda e^t -1} = \sum_{n=0}^\infty B_n (\lambda)
\dfrac{t^n}{n!} . & (3a)}
$$
If $\lambda =e^{2 \pi i /f}$, $f\in \Bbb N$ and $\lambda \in\Bbb
C$, then Eq-(3) is reduced to (3a). Eq-(3a) is obtained by Kim
[3]. Let $u\in\Bbb C$, then by substituting $x=0$ into Eq-(A1),
we set
$$\eqalignno{ &
\dfrac{1-u}{e^t -u} =\sum_{n=0}^\infty H_n (u)
\dfrac{t^n}{n!}.&(3b)}
$$
$H_n (u)$ is denoted Frobenius-Euler numbers. Relation between
$H_n (u)$ and $B_n (\lambda)$  is given by the following theorem:

\proclaim{Theorem 1} Let $\lambda\in\Bbb Z_p$. Then
$$\eqalignno{ &
B_n (\lambda) =\dfrac{ \log \lambda}{\lambda -1} H_n
(\lambda^{-1}) + \dfrac{n H_{n-1}(\lambda^{-1} )}{ \lambda
-1},&(4)\cr & B_0 (\lambda) =\dfrac{\log \lambda }{\lambda -1} H_0
(\lambda^{-1}).}
$$
\endproclaim
\demo{Proof} By using Eq-(3), we have
$$
\split & \sum_{n=0}^\infty B_n (\lambda ) \dfrac{t^n}{n!}
=\dfrac{\log \lambda +t}{\lambda e^t -1 } = \dfrac{\log \lambda
}{\lambda e^t -1} + \dfrac{t}{\lambda e^t -1}\cr &=
\dfrac{1-\lambda^{-1}}{(1-\lambda^{-1})\lambda}\cdot \left(
 \dfrac{\log \lambda }{e^t -\lambda^{-1}}  \right) -
 \dfrac{(1-\lambda^{-1})}{(e^t -\lambda^{-1})}\cdot \dfrac{t}{\lambda
 (1-\lambda^{-1})}\cr
 & =\dfrac{\log \lambda}{\lambda -1} \sum_{n=0}^\infty H_n
 (\lambda^{-1}) \dfrac{t^n}{n!} + \dfrac{t}{\lambda -1} \sum_{n=0}^\infty
 H_n (\lambda^{-1})\dfrac{t^n}{n!} ,
\endsplit
$$
the next  to the last step being  a consequence of Eq-(3b).  After
some  elementary calculations, we have
$$
\sum_{n=0}^\infty B_n (\lambda ) \dfrac{t^n}{n!} =\dfrac{\log
\lambda}{\lambda -1} H_0 (\lambda^{-1}) + \sum_{n=1}^\infty \left(
\dfrac{\log \lambda}{\lambda -1}H_n (\lambda^{-1}) +
\dfrac{n}{\lambda -1} H_{n-1}
(\lambda^{-1})\right)\dfrac{t^n}{n!}.
$$
By comparing  coefficient $\dfrac{t^n}{n!}$ in the above, then we
obtain the desired result. \quad\qed
\enddemo

Observe that, if $\lambda \in T_p$ in Eq-(4), then we have,
$B_0(\lambda )=0$ and $B_n (\lambda )= \dfrac{n H_{n-1}
(\lambda^{-1})}{\lambda -1}$,  $n\geq 1$.

By  Eq-(3) and Eq-(4), we obtain the following formula:

 For $n\geq 0$, $\lambda \in \Bbb Z_p$

 $$
\eqalignno{ & \int_{\Bbb Z_p} \lambda^x x^n d\mu_{1} (x) = \cases
\dfrac{\log \lambda}{\lambda -1} H_0 (\lambda^{-1}), & \text{   $n=0$}\\
\dfrac{\log \lambda}{\lambda -1} H_n (\lambda^{-1})
+\dfrac{n}{\lambda -1} H_{n-1} (\lambda^{-1}), & \text{ $n>0$}
\endcases &(4a)}
 $$

and

$$
\eqalignno{ & \int_{\Bbb Z_p} \lambda^x x^n d\mu_{1} (x) =B_n
(\lambda), \quad n \geq 0. &(4b)}
 $$
Now, we define $\lambda -$Bernoulli polynomials, we use these
polynomials to give the sums powers of consecutive. The
$\lambda$-Bernoulli polynomials are defined by means of the
following generating function:

$$\eqalignno{ &\dfrac{\log \lambda +t}{ \lambda e^t -1}
 e^{tx}=\sum_{n=0}^\infty B_n (\lambda ; x) \dfrac{t^n}{n!} . &(5)}
$$
By Eq-(3) and Eq-(5), we have
$$
B_n (\lambda ; x)= \sum_{k=0}^n \binom{n}{k} B_k
(\lambda)x^{n-k}.$$

The Witt's formula for $B_n (\lambda; x)$ is given by the
following theorem:

\proclaim{Theorem 2} For $k \in \Bbb N$ and $\lambda \in \Bbb
Z_p$, we have
$$\eqalignno{ &B_n (\lambda; x) =\int_{\Bbb Z_p} (x+y)^n \lambda^y d\mu_1 (y) . &(6)}
$$
\endproclaim
\demo{Proof} By substituting $f(y) = e^{t(x+y)}\lambda^y$ into
Eq-(1), we have
$$
\int_{\Bbb Z_p} e^{t(x+y)}\lambda^y d\mu_1 (y) =\sum_{n=0}^\infty
B_n (\lambda; x) \dfrac{t^n}{n!} =\dfrac{(\log\lambda
+t)e^{tx}}{\lambda e^t -1}.
$$

By using Taylor expansion of $e^{tx}$ in the left side of the
above equation, after some elementary calculations, we obtain the
desired result. \quad\qed
\enddemo

We now give the distribution of the $\lambda-$Bernoulli
polynomials.

\proclaim{Theorem 3} Let $n\geq 0$, and let $d\in\Bbb Z^+$. Then
we have
$$\eqalignno{ & B_n (\lambda; x) =d^{n-1} \sum_{a=0}^{d-1} \lambda^a B_n \left(\lambda^d ; \dfrac{x+a}{d} \right) . &(7)}
$$
\endproclaim
\demo{Proof} By using  Eq-(6),
$$\split
B_n (x; \lambda) &=\int_{\Bbb Z_p} (x+y)^n \lambda^y d \mu_1
(y)\cr &= \lim_{N \to \infty } \dfrac{1}{dp^N} \sum_{y=0}^{dp^N
-1} (x+y)^n \lambda^y\cr &= \lim_{N \to \infty } \dfrac{1}{dp^N}
\sum_{a=0}^{d -1} \sum_{y=0}^{p^N -1} (a+dy+x)^n \lambda^{a+dy}\cr
&= d^{n-1} \lim_{N \to \infty }  \dfrac{1}{p^N} \sum_{a=0}^{d -1}
\lambda^a \sum_{y=0}^{p^N -1} \left(\dfrac{a+x}{d} +y\right)^n
(\lambda^{d})^y \cr &= d^{n-1}  \dfrac{1}{p^N} \sum_{a=0}^{d -1}
\lambda^a \int_{\Bbb Z_p} \left(\dfrac{a+x}{d} +y\right)^n
(\lambda^{d})^y .
\endsplit
$$
Thus, we have the desired result.\quad\qed
\enddemo

By substituting $x=0$ into Eq-(7), we have the following
corollary:

\proclaim{Corollary 1} For $m,n\in \Bbb N$, we have
$$\eqalignno{ & m B_n (\lambda) = \sum_{j=0}^{n}\binom{n}{j} B_j (\lambda^m ) m^j \sum_{a=0}^{m-1}
\lambda^a a^{n-j} . &(8)}
$$
( Gauss multiplicative formula for $\lambda$-Bernoulli numbers).
\endproclaim
By Eq-(8), we have

\proclaim{Theorem  4} For $m,n\in \Bbb N$ and $\lambda \in\Bbb
Z_p$, we have
$$\eqalignno{ & m B_n (\lambda) -m^n [m]_\lambda B_n (\lambda^m )=
 \sum_{j=0}^{n-1}\binom{n}{j} B_j (\lambda^m ) m^j \sum_{k=1}^{m-1}
\lambda^k k^{n-j} . &(9)}
$$
\endproclaim

\proclaim{Theorem  5} Let $k\in\Bbb Z$, with $k>1$. Then we have
$$\eqalignno{ &  B_l (\lambda;k) -\lambda^{-k} B_l(\lambda) =\lambda^{-k}
 l\sum_{n=0}^{k-1}\lambda^n n^{l-1} +
(\lambda^{-k} \log\lambda )  \sum_{n=0}^{k-1} n^l \lambda^l .
&(10)}
$$
\endproclaim

\demo{Proof} We set
$$\eqalignno{ & -\sum_{n=0}^\infty e^{(n+k)t}\lambda^n
+\sum_{n=0}^\infty e^{nt}\lambda^{n-k} =\sum_{n=0}^\infty
e^{nt}\lambda^{n-k}&(10a)\cr &= \sum_{l=0}^\infty (\lambda^{-k}
\sum_{n=0}^{k-1} n^{l} \lambda^n )\dfrac{t^l}{l!}\cr &=
\sum_{l=1}^\infty (\lambda^{-k}l \sum_{n=0}^{k-1} n^{l-1}
\lambda^n )\dfrac{t^{l-1} }{l!}.}
$$
Multiplying $(t +\log \lambda)$ both side of Eq-(10a), then by
using Eq-(3) and Eq-(5), after some elementary calculations, we
have
$$\eqalignno{ & \sum_{l=0}^\infty ( B_l (\lambda;k) -\lambda^{-k}  B_l (\lambda)) \dfrac{t^l}{l!}\cr
&= \sum_{l=0}^\infty (\lambda^{-k} l  \sum_{n=0}^{k-1}\lambda^n
n^{l-1} +\lambda^{-k}\log \lambda  \sum_{n=0}^{k-1}n^l \lambda^l )
\dfrac{t^l}{l!} . &(10b)}
$$
By comparing coefficient $\dfrac{t^l}{l!} $ in both sides of
Eq-(10b). Thus we arrive at the Eq-(10). Thus we complete the
proof of theorem.\quad\qed

\enddemo

Observe that $\lim_{\lambda \to 1} B_l (\lambda) =B_l$. For
$\lambda \to 1$, then Eq-(10) reduces the following:

$$
B_l (k) -B_l = l \sum_{n=0}^{k-1} n^{l-1}.
$$
If $\lambda \in T_p$, then Eq-(10) reduces to the following
formula:

$$
B_l (\lambda; k) -\lambda^{-k} B_l (\lambda) = \lambda^{-k} l
\sum_{n=0}^{k-1} \lambda^n  n^{l-1}.
$$

\remark{Remark} Garrett and Hummel [1B] proved combinatorial proof
of $q$-analogue of $$\sum_{k=1}^n k^3= {\binom{n+1}{k}}^2$$ as
follows:
$$
\sum_{k=1}^n q^{k-1} [k]_q^2 \left( {{k-1} \brack {2}}_{q^2} +
{{k+1} \brack {2}}_{q^2} \right) = {{n+1} \brack {2}}_{q}^2 ,
$$
where ${{n} \brack {k}}_{q}= \prod_{j=1}^k \dfrac{[n+1-j]_q
}{[j]_q}$, $q$-binomial coefficients. Garrett and Hummel, in their
paper, asked for a simpler $q$-analogue of the sums of cubes. As a
response to Garrett and Hummel's question, in [11], Kim
constructed the following formula
$$\split
S_{n,q^h} (k) & =\sum_{l=0}^{k-1} q^{h^l} [l]^n \cr
  &= \dfrac{1}{n+1} \sum_{j=0}^n \binom{n+1}{j}
  \beta_{j,q} q^{kj} [k]^{n+1-j}  - \dfrac{(1-q^{(n+1)k})\beta_{n+1,q}}{n+1}, \endsplit
$$
where $\beta_{j,q}$ are the $q$-Bernoulli numbers which were
defined by
$$
e^{\frac{t}{1-q}} \dfrac{q-1}{\log q} -t \sum_{n=0}^\infty q^{n+x}
e^{[n+x]t} = \sum_{n=0}^\infty \dfrac{\beta_{n,q} (x)}{n!} t^n ,
\quad \vert q\vert <1, \vert t \vert <1,
$$
$\beta_{n,q} (0) =\beta_{n,q}$. cf. [11]

Schlosser [20] gave for $n=1,2,3,4,5$ the value of $S_{n,q^h}[k]$.
In [28], the authors also gave  another proof of $S_{n,q}(k)$
formula.
\endremark

{\bf\centerline {\S 2. Hurwitz's type $\lambda$-zeta function}}

 \vskip 20pt

In this section, by using generating function of
$\lambda$-Bernoulli polynomials,  we construct Hurwitz's type
$\lambda$-zeta function, which is interpolate $\lambda$-Bernoulli
polynomials at negative integers. By Eq-(5), we get
$$\split
F_\lambda (t;x ) &= \dfrac{\log \lambda +t}{\lambda e^t -1} e^{xt}
= -(\log \lambda +t )\sum_{n=0}^\infty \lambda^n e^{(n+x)t}\cr &=
\sum_{n=0}^\infty B_n (\lambda ) \dfrac{t^n}{n!}. \endsplit
$$

By using $\dfrac{d^k}{dt^k}$ derivative operator to the above, we
have
$$\split
B_k (\lambda; x) &= \left. \dfrac{d^k}{dt^k} F_\lambda
(t;x)\right\vert_{t=0},\cr B_k (\lambda; x) &=
 -\log \lambda \sum_{n=0}^\infty \lambda^n (n+x)^k -k
 \sum_{n=0}^\infty (n+x)^{k-1} \lambda^n .
\endsplit
$$
Thus we arrive at the following theorem:

\proclaim{Theorem 6} For $k\geq 0$, we have
$$
- \dfrac{1}{k} B_k (\lambda ;x) =\dfrac{\log \lambda^k}{k}
\sum_{n=0}^\infty \lambda^n (n+x)^k + \sum_{n=0}^\infty \lambda^n
(n+x)^{k-1} .
$$
\endproclaim

Consequently, we define Hurwitz type zeta function as follows:

\definition{Definition 1} Let $s\in \Bbb C$. Then we define
$$\eqalignno{ &
\zeta_\lambda (s, x) = \dfrac{\log \lambda}{1-s}
 \sum_{n=0}^\infty \dfrac{ \lambda^n }{(n+x)^{s-1}} +
 \sum_{n=0}^\infty \dfrac{ \lambda^n }{(n+x)^{s}}.  &(11)}
$$

\enddefinition

Note that $\zeta_\lambda (s, x)$ is analytic continuation, except
for $s=1$, in whole complex plane. By Definition 1 and Theorem 6,
we have the following:

 \proclaim{Theorem 7} Let $s=1-k$, $k \in \Bbb N $,
$$\eqalignno{ &
\zeta_\lambda (1-k , x) =- \dfrac{B_k (\lambda, x)}{ k}. &(12)}
$$
\endproclaim

\vskip 20pt

{\bf\centerline {\S 3. Generalized $\lambda$-Bernoulli numbers}}

{\bf\centerline { associated with Dirichlet type
$\lambda$-$L$-functions}}

\vskip 20pt

By using Eq-(0), we define
$$\eqalignno{ &
I_1 (f_d ) =I_1 (f) + \sum_{j=0}^{d-1} f^\prime (j),&(12)}
$$
where $f_d (x) = f(x+d)$, $\int_{\Bbb X} f(x) d\mu (x) = I_1 (f).$

Let $\chi$ be Dirichlet character with conductor $d\in \Bbb
N^{+}$, $\lambda \in \Bbb Z_p$.

By substituting $f(x) =\lambda^x \chi (x) e^{tx}$ into Eq-(12),
then we have
$$\eqalignno{
  \int_{\Bbb X} \chi (x) \lambda^{x} e^{tx} d \mu_1 (x) &=
\sum_{j=0}^{d-1} \dfrac{\chi (j) \lambda^j e^{tj } (\log \lambda
+t) }{\lambda^d e^{dt} -1}\cr &=\sum_{n=0}^{\infty} B_{n,\chi}
(\lambda )\dfrac{t^n}{n!} .&{(12a)}}
$$

By Eq-(12a),  we easily see that
$$\eqalignno{ &
B_{n,\chi} (\lambda ) = \int_{\Bbb X} \chi (x) x^n  \lambda^x d
\mu_1 (x). &{(12b)}}
$$

From Eq-(12a), we define generating function of generalized
Bernoulli number by
$$\eqalignno{ &
F_{\lambda , \chi} (t) =
 \sum_{j=0}^{d-1} \dfrac{ \chi (j) \lambda^j e^{tj}(\log \lambda +t) }{\lambda^d e^{dt} -1}
 =\sum_{n=0}^{\infty} B_{n} (\lambda) \dfrac{t^n}{n!}. &(12c)}
$$
Observe that if $\lambda \in T_p$, then the above formula reduces
to
$$
F_{\lambda , \chi} (t) =\sum_{j=0}^{d-1} \dfrac{ \chi (j)
\lambda^j e^{tj}t }{\lambda^d e^{dt} -1} =\sum_{j=0}^{\infty}
B_{n} (\lambda) \dfrac{t^n}{n!}\quad\text{( for detail see cf.
[3,16,18,22,23,24]).}
$$
From the above, we easily see that
$$
F_{\lambda , \chi} (t) =- (\log \lambda +t) \sum_{m=1}^\infty
\chi(m)\lambda^m e^{tm} = \sum_{n=0}^\infty B_{n,\chi} (\lambda)
\dfrac{t^n}{n!}.
$$
By applying $\left. \dfrac{d^k}{dt^k}\right\vert_{t=0}$ derivative
operator both sides of the above equation, we arrive at the
following theorem:

\proclaim{Theorem 8} Let $k\in \Bbb Z^{+}$, $\lambda \in\Bbb Z_p$
and let $\chi$ be Derichlet character with conductor $d$. Then we
have
$$\eqalignno{ &
\sum_{m=1}^\infty \chi (m) \lambda^m m^{k-1} + \dfrac{\log
\lambda}{ k} \sum_{m=1}^\infty \lambda^m \chi(m) m^k = -
\dfrac{B_{k,\chi}(\lambda)}{k}. &(13)}
$$
\endproclaim

\definition{Definition 2( Dirichlet type $\lambda$-$L$ function)} For
$\lambda, s\in\Bbb C$, we define
$$\eqalignno{ &
L_{\lambda } (s,\chi) =\sum_{m=1}^\infty \dfrac{\lambda^m
\chi(m)}{ m^s} - \dfrac{ \log\lambda}{s-1} \sum_{m=1}^\infty
\dfrac{\lambda^m \chi(m)}{ m^{s-1}}.&(14)}
$$

\enddefinition

Relation between $L_\lambda (s,\chi )$  and $\zeta_\lambda (s,y)$
is given by the following theorem :

\proclaim{Theorem 9} Let  $s\in\Bbb C$ and $d\in \Bbb Z^+$. Then
we have
$$\eqalignno{ & L_\lambda (s,\chi ) =d^{-s}
\sum_{a=1}^d \lambda^a \chi(a) \zeta_{\lambda^d} \left(s,
\dfrac{a}{d} \right).}
$$
\endproclaim

\demo{Proof} By substituting $m=a +dk$,  $a=1,2,\cdots, d$,
$k=0,1,\cdots, \infty$, into Eq-(14), we have
$$
\split &
 L_\lambda (s,\chi ) = \sum_{a=1}^d \sum_{k=0}^\infty
 \dfrac{\lambda^{a+dk} \chi (a+dk )}{(a+dk)^s } -\dfrac{\log \lambda }{s-1}
\sum_{a=1}^d \sum_{k=0}^\infty  \dfrac{\lambda^{a+dk} \chi (a+dk
)}{(a+dk)^{s-1} }\cr &\quad = d^{-s} \sum_{a=1}^d (\lambda^a \chi
(a)  ) \left[ \sum_{k=0}^\infty
\dfrac{(\lambda^d)^k}{(k+\frac{a}{d})^s} - \dfrac{\log \lambda^d
}{s-1} \sum_{k=0}^\infty
\dfrac{(\lambda^d)^k}{(k+\frac{a}{d})^{s-1} }\right].
\endsplit
$$

By  using Eq-(11) in the above  we obtain the desired result.
\quad\qed
\enddemo

\proclaim{Theorem 10} For $k\in\Bbb Z^+$,  we have
$$\eqalignno{ & L_\lambda (1-k,\chi ) =- \dfrac{1}{k} B_{k,\chi} (\lambda), \ \ k>0.}
$$
\endproclaim

\demo{Proof}  By substituting $s=1-k$ in Definition 2 and using
Eq-(13), we easily obtain the desired result.$\quad\qed$
\enddemo

\remark{Remark} If $\lambda\in T_p$, then from Definition 2, we
have
$$
L_\lambda( s,\chi) =\sum_{m=1}^\infty \dfrac{\lambda^m
\chi(m)}{m^s}.
$$

In [18], Koblitz studied on this function. He gave the name of
this function ``twisted $L$-function"  for $\lambda$ is $r$-th
root of $1$.  In [22,23,24], Simsek studied of this functions. He
gave fundamental properties of this function as well.

In [16], Kim et all gave $\lambda$-$(h,q)$ zeta function and
$\lambda$-$(h,q)$ $L$-function. These functions interpolate
$\lambda-(h,q)$-Bernoulli numbers. Observe that, if we take
$s=1-k$ in Theorem 9, and then using Eq-(12) in Theorem 7, we get
another proof of Theorem 10.
\endremark

\vskip 20pt

{\bf\centerline {\S 4. $\lambda$-Bernoulli numbers of order $r$
associated with multiple zeta function}}
\vskip 20pt

In this section, we define generating function of
 $\lambda$-Bernoulli numbers of order $r$. By using Mellin
 transforms and Cauchy residue theorem, we obtain  multiple
 zeta function which is given in Eq-(C). We also gave relations
 between $\lambda$-Bernoulli polynomials  of order $r$ and multiple
 zeta function at negative integers. This relation is important
 and very interesting.  Let $r\in\Bbb Z^+$. Generating function of
 $\lambda$-Bernoulli numbers of order $r$ is defined by
 $$
\eqalignno{ & F_\lambda^{(r)} (t) = \left( \dfrac{\log \lambda
+t}{\lambda e^t -1} \right)^r = \sum_{n=0}^\infty  B_n^{(r)}
(\lambda) \dfrac{t^n}{n!}.&(15)}
 $$
Generating function of $\lambda$-Bernoulli polynomials of order
$r$ is defined by
$$
F_\lambda^{(r)} (t,x) = F_\lambda^{(r)} (t) e^{tx} =
\sum_{n=0}^\infty B_n^{(r)} (\lambda) \dfrac{t^n}{n!}.
$$

Observe that when $r=1$, Eq-(15) reduces to Eq-(3). By applying
Mellin transforms to the Eq-(15) we get
$$
\dfrac{1}{\Gamma (s)} \int_0^\infty \lambda^r e^{-tr}
F_\lambda^{(r)} (-t) (t-\log \lambda )^{s-r-1} dt= \sum_{n_1
,\cdots, n_r =0}^\infty \dfrac{1}{(n_1 + n_2 +\cdots + n_r +r
)^s}.
$$
Thus, we get, by (C)
$$
\zeta_r (s) =\dfrac{1}{\Gamma (s)} \int_0^\infty \lambda^r e^{-tr}
F_\lambda^{(r)} (-t) (t-\log \lambda )^{s-r-1} dt .
$$

By using the above relation, we obtain the following theorem:

\proclaim{Theorem 11} Let $r,m\in \Bbb Z^+$.  Then we have
$$\eqalignno{ &
\zeta_r (-m) =(-\lambda)^r m! \sum_{j=0}^\infty \binom{-m-r-1}{j}
(\log \lambda )^j \dfrac{B_{m+r+j}(\lambda; r)}{(m+r+j)!}. &(D1)}
$$
\endproclaim

\remark{Remark} If $\lambda \to 1$, the above theorem reduces to
$$\eqalignno{ &
\zeta_r (-m) =(-1)^r m! \dfrac{B_{m+r}(1; r)}{(m+r)!} &(D2)}
$$
which is given Theorem 6 in [12].

By (D1) and (D2), we obtain  relation between $\lambda$-Bernoulli
polynomials of order $r$ and ordinary Bernoulli  polynomials of
order $r$ as follows:

$$
B_{m+r} (r) =\lambda^r \sum_{j=0}^\infty \binom{-m-r-1}{j} (\log
\lambda)^j \dfrac{B_{m+r+j}(\lambda; r)}{(m+r+j)!}(m+r)!
$$
where $m,r \in \Bbb Z^+$.
\endremark

We now give relations between $B_n^{(r)}(\lambda)$ and $H_n^{(r)}
(\lambda^{-1})$ as follows:

If $\lambda \in T_p$, then Eq-(15) reduces to the following
equation
$$
\dfrac{t^r}{ (\lambda e^t -1)^r} =\sum_{n=0}^\infty
B_n^{(r)}(\lambda) \dfrac{t^n}{n!}.
$$
Thus by the above equation, we easily see that
$$\split
t^r &= (\lambda e^t -1)^r e^{B^{(r)}(\lambda)t} \cr
  &= \sum_{l=0}^r \lambda^l (-1)^{r-l} e^{(B^{(r)}(\lambda) +l)t
  }\cr
  &=\sum_{n=0}^\infty ( \sum_{l=0}^r \lambda^l (-1)^{r-l} (
  B^{(r)} (\lambda ) +l)^n  )\dfrac{t^n}{n!}.
\endsplit
$$

Consequently we have
$$
 \sum_{l=0}^r \lambda^l (-1)^{r-l} (B^{(r)} (\lambda )+l)^n =
 \cases
0 & \text{if\   $n\neq r$}\\
1 & \text{if \  $n=r$ .}
\endcases
$$

By Eq-(15) we obtain
$$
\sum_{n=0}^\infty B_n^{(r)} (\lambda ) \dfrac{t^n}{n!} =
\dfrac{t^r}{ (\lambda -1)^r} \sum_{n=0}^\infty H_n^{(r)}
(\lambda^{-1}) \dfrac{t^n}{n!}.
$$
By comparing coefficient $\dfrac{t^n}{n!}$ in the both sides of
the above equation, we have
$$
 B_{n+r}^{(r)} (\lambda ) = \dfrac{\Gamma (n+r +1)}{\Gamma (n +1)}
 \dfrac{1}{(\lambda-1)^r} H_n^{(r)} (\lambda^{-1}).
$$
Observe that, if we take $r=1$, then the above identity reduce to
Eq-(4), that is
$$
B_{n+1} (\lambda) =\dfrac{(n+1)}{\lambda -1} H_n (\lambda^{-1}).
$$

\vskip 20pt

{\bf\centerline {\S 5. $\lambda$-Bernoulli numbers and polynomials
associated}}

{\bf\centerline { with multivariate $p$-adic invariant integral}}
\vskip 20pt

In this section, we give generalized $\lambda$-Bernoulli numbers
of order $r$. Consider the multivariate $p$-adic invariant
integral
 on
$\Bbb Z_p$ to define $\lambda$-Bernoulli numbers and polynomials.
$$\eqalignno{ &
\underbrace{\int_{\Bbb Z_p}\cdots\int_{\Bbb
Z_p}}_{r-\text{times}}\lambda^{w_1 x_1 +\cdots+w_r x_r} e^{(w_1
x_1 +\cdots+w_r x_r)t} d\mu_1 (x_1 ) \cdots d\mu_1 (x_r)\cr
&=\dfrac{ ( w_1 \log \lambda +w_1 t)\cdots (w_r \log \lambda +w_r
t) }{(\lambda^{w_1} e^{w_1 t} -1 )\cdots (\lambda^{w_r} e^{w_r t}
-1 ) }&(16)\cr & = \sum_{n=0}^\infty B_n^{(r)} (\lambda; w_1 , w_2
, \cdots, w_r ) \dfrac{t^n}{n!},}
$$
where we called $B_n^{(r)} (\lambda; w_1 , w_2 , \cdots, w_r ) $
$\lambda$-extension of Bernoulli numbers. Substituting $\lambda
=1$ into Eq-(16), $\lambda$-extension of Bernoulli numbers  reduce
to Barnes Bernoulli numbers as follows :

$$\dfrac{ ( w_1 t)\cdots (w_r
t) }{( e^{w_1 t} -1 )\cdots ( e^{w_r t} -1 ) } = \sum_{n=0}^\infty
B_n^{(r)} (w_1 ,\cdots, w_r )\dfrac{t^n}{n!},
$$
where  $B_n^{(r)} (w_1 ,\cdots, w_r )$ are denoted Barnes
Bernoulli umbers and $w_1, \cdots, w_r$ complex numbers with
positive real parts [1A,7,27]. Observe that when $w_1 =w_2 =\cdots
= w_r =1$ in Eq-(16), we obtain the $\lambda$-Bernoulli numbers of
higher order as follows:

$$
\left( \dfrac{\log \lambda +t }{\lambda e^t -1} \right)^r
=\sum_{n=0}^\infty B_n^{(r)} (\lambda) \dfrac{t^n}{n!}.
$$
We note that $B_n^{(r)} (\lambda; 1,1,\cdots, 1)=B_n^{(r)}
(\lambda)$.

Consider
$$
\left( \dfrac{\log \lambda +t }{\lambda e^t -1} \right)^r e^{xt}
=\sum_{n=0}^\infty B_n^{(r)} (\lambda; x) \dfrac{t^n}{n!}.
$$

Observe that
$$
\split & \sum_{n=0}^\infty B_n^{(r)} (\lambda; x) \dfrac{t^n}{n!}
=\left( \dfrac{\log \lambda +t }{\lambda e^t -1} \right)^r
e^{(\log \lambda +t )x} \lambda^{-x}\cr
 & = \dfrac{1}{\lambda^x} \sum_{m=0}^\infty B_m^{(r)} (\lambda; x)
 \dfrac{(t+\log \lambda)^m}{m!}\cr
 &= \dfrac{1}{\lambda^x} \sum_{m=0}^\infty \dfrac{B_m^{(r)} (\lambda; x) }{m!}
\sum_{l=0}^m \binom{m}{l} (\log \lambda)^m t^{m-l}\cr &=
\sum_{m=0}^\infty \left(
 \dfrac{1}{\lambda^x} \sum_{l=0}^\infty \dfrac{B_{n+l}^{(r)}
 (\lambda;x)}{l!}(\log \lambda)^l
 \right) \dfrac{t^n}{n!}.
\endsplit
$$
Now, comparing coefficient $\dfrac{t^n}{n!}$ both sides of the
above equation, we easily arrive at the following theorem:

\proclaim{Theorem 12} For $n,r\in \Bbb N$ and $\lambda \in\Bbb
Z_p$, we have
$$
B_n^{(r)} (\lambda;x) = \dfrac{1}{\lambda^r } \sum_{l=0}^\infty
B_{n+l}^{(r)} (\lambda; x) \dfrac{(\log \lambda)^l}{l!},
$$
where $0^l =\cases
1 & \text{if\   $l=0$}\\
0 & \text{if \  $l\neq 0$ .}
\endcases
$
\endproclaim

\remark{Remark}
 In Theorem 12, we see that
 $$
\lim_{\lambda \to 1} B_n^{(r)} (\lambda ;x) =\cases
B_n^{(r)} (x) & \text{if\   $l=0$,}\\
0 & \text{if \  $l\neq 0$. }
\endcases
 $$
\endremark

\vskip 20pt

{\bf\centerline {\S 6. $\lambda$-Bernoulli numbers and polynomials
in}}

{\bf\centerline { the  space of locally constant}} \vskip 20pt

In this section, we construct partial $\lambda$-zeta functions, we
need this function in the following section. We need this function
in the following section.
 By Eq-(3b),
Frobenius-Euler polynomials are defined by means of the following
generating function:
$$
\left( \dfrac{1-u}{e^t -u}
 \right) e^{xt} =\sum_{n=0}^\infty H_n (u,x) \dfrac{t^n}{n!}.
$$
As well known, we note that the Frobenius-Euler polynomials of
order  $r$ were defined by
$$\left( \dfrac{1-u}{e^t -u}
 \right)^r e^{xt} =\sum_{n=0}^\infty H_n^{(r)} (u,x) \dfrac{t^n}{n!}.
$$
The case $x=0$, $H_n^{(r)} (u,0) =H_n^{(r)} (u)$, which are called
Frobenius-Euler numbers of order $r$.

If $\lambda \in T_p$, then $\lambda$-Bernoulli polynomials of
order $r$ are given by
$$
\dfrac{t^r}{ (\lambda e^t -1)^r} e^{xt} =  \sum_{n=0}^\infty
B_n^{(r)} (\lambda ;x) \dfrac{t^n}{n!}.
$$
Hurwitz type  $\lambda$-zeta function is given by
$$\eqalignno{ &
\zeta_\lambda (s,x) = \sum_{n=0}^\infty \dfrac{\lambda^n
}{(n+x)^s} , \quad
 \lambda \in T_p . &(17)}
$$

Thus, from Theorem 7, we have
$$\eqalignno{ &
\zeta_\lambda  (1-k,x) = -\dfrac{1}{k}B (\lambda ;x),\quad \quad
 k \in \Bbb Z^+ . &(17a)}
$$
We now define $\lambda$-partial zeta function as follows
$$\eqalignno{ &
H_\lambda (s,a\vert F)= \sum_{m\equiv a \pmod{F}}
\dfrac{\lambda^m}{m^s}.&(17b)}
$$
From (17), we have
$$\eqalignno{ &
H_\lambda (s,a\vert F)= \dfrac{\lambda^a }{F^s}
\zeta_{\lambda^F}\left(s, \dfrac{a}{F}\right),&(17c)}
$$
where $\zeta_{\lambda^F}\left(s, \dfrac{a}{F}\right)$ is given by
Eq-(17). By Eq-(17a) we have
$$\eqalignno{ &  H_\lambda (1-n ,a\vert F)
 = -\dfrac{F^{n-1} \lambda^a B_n (\lambda^F ; \dfrac{a}{F})}{n},\quad \quad
 n \in \Bbb Z^+ . &(18)}
$$

If $\lambda \in T_p $, then by Eq-(14), we have

$$
L_\lambda (s,\chi ) =\sum_{n=1}^\infty \dfrac{\lambda^n \chi
(n)}{n^s} ,
$$
where $s\in\Bbb C$, $\chi$ be the primitive Dirichlet character
with conductor $f\in\Bbb Z^+$. By Theorem 9, Eq-(17c) and  Eq-(18)
we easily see that
$$
L_\lambda (s,\chi ) =\sum_{a=1}^F  \chi (a) H_\lambda \left(s,
\dfrac{a}{F}\right),
$$
and
$$L_\lambda (1-k ,\chi ) = -\dfrac{B_{k,\chi}(\lambda)}{k}, \quad
k\in \Bbb Z^+ ,
$$
where $B_{k,\chi}(\lambda)$ is defined by

$$
\sum_{a=0}^{F-1} \dfrac{t \lambda^a \chi (a) e^{at}}{ \lambda^F
e^{Ft} -1 } =\sum_{a=0}^{\infty} B_{n,\chi} (\lambda)
\dfrac{t^n}{n!}, \quad \lambda \in T_p
$$
and $F$ is multiple of $f$.

\remark{Remark}
$$
\dfrac{B_m (\lambda)}{m} =\dfrac{1}{\lambda-1} H_{n-1}
(\lambda^{-1}), \quad \lambda\in T_p .$$
\endremark

\vskip 20pt

{\bf\centerline {\S 7. $p$-adic interpolation function}}
 \vskip 20pt

In this section we give $p$-adic $\lambda$-$L$ function. Let $w$
be the Teichimuller character and let $<x> =\dfrac{x}{w(x)}$.

When $F$ is multiple of $p$ and $f$ and $(a,p)=1$, we define
$$
H_{p,\lambda} (s, a\vert F) = \dfrac{1}{s-1} \lambda^a <a>^{1-s}
\sum_{j=0}^\infty \binom{1-s}{j} \left( \dfrac{F}{a} \right)^j B_j
(\lambda^F ).
$$

From this we note that
$$\split
H_{p,\lambda} (1-n, a\vert F) &= -\dfrac{1}{n}
\dfrac{\lambda^a}{F} <a>^n \sum_{j=0}^n \binom{n}{j} \left(
\dfrac{F}{a} \right)^j B_j (\lambda^F )\cr &= - \dfrac{1}{n}F^{n-1
} \lambda^a w^{-n} (a) B_n (\lambda^F ; \dfrac{a}{F})\cr &=w^{-n}
(a) H_\lambda (1-n ; \dfrac{a}{F}),\endsplit$$

since by Theorem 3 for $\lambda \in T_p$, Eq-(18).

By using this formula, we can consider $p$-adic
$\lambda$-$L$-function for $\lambda$-Bernoulli numbers as follows:
$$
L_{p,\lambda} (s,\chi) =\sum_{a=1 \atop p\not| a}^F \chi (a)
H_{p,\lambda} \left( s, \dfrac{a}{F} \right) .
$$

By using the above definition, we have
$$\split
L_{p,\lambda} (1-n, \chi) &=\sum_{a=1 \atop p\not\vert a}^F \chi
(a)H_{p,\lambda} \left( 1-n, \dfrac{a}{F} \right)\cr
&=
-\dfrac{1}{n} \left(
 B_{n,\chi w^{-n} } (\lambda) -p^{n-1} \chi w^{-n} (p)
 B_{n,\chi w^{-n}} (\lambda^p)
\right).
\endsplit
$$
Thus, we define the formula

$$
L_{p,\lambda} (s,\chi) =\dfrac{1}{F} \dfrac{1}{s-1} \sum_{a=1}^F
\chi (a) \lambda^a <a>^{1-s} \sum_{j=0}^\infty \binom{1-s}{j}B_j
(\lambda^F)
$$
for $s\in \Bbb Z_p$.

\vskip 20pt

 \vskip 10pt \quad  Yilmaz  Simsek

 \quad  University of Akdeniz, Faculty of  Art and Science,

 \quad Department of Mathematics 07058 Antalya, Turkey

 \quad e-mail: ysimsek$\@$akdeniz.edu.tr \vskip 10pt

 \vskip 10pt
\quad Taekyun Kim

\quad Jangjeon Research Institute for Mathematical Sciences and
Physics,

\quad 252-5 Hapcheon-Dong Hapcheon-Gun Kyungnam, 678-802, S. Korea

\quad e-mail:\text{tkim64$\@$hanmail.net, tkim$\@$kongju.ac.kr}

\vskip 10pt

\quad    Daeyeoul Kim

\quad    Department of Mathematics and Institute of

\quad    Pure and Applied Mathematics,

\quad     Chonbuk National Univ., Chonju, 561-756, Korea

\quad   e-mail: daeyeoul\@chonbuk.ac.kr \vskip 20pt


\Refs \widestnumber\key{999999}




\ref \key 1A \by  E\. W\. Barnes
  \paper On the theory of the multiple gamma functions
 \jour  Trans. Camb. Philos. Soc.
 \yr 1904
\pages 374-425 \vol 19 \endref


\ref \key 1B \by  K\. C\. Garret and K\. Hummel
  \paper A combinatorial proof of the sum of $q$-cubes
 \jour  Electro. J. Comb. 11\# R
 \yr 2004
\pages  \vol 9 \endref



\ref \key 2
 \by  K. Iwasawa
 \book  Lectures on $p$-adic $L$-function
\publ Princeton Univ. \yr 1972
\endref



\ref \key 3
 \by  T. Kim
  \paper  An analogue of Bernoulli numbers and their congruences
 \jour  Rep. Fac. Sci. Engrg. Saga Univ. Math.
 \yr 1994
\pages 21--26 \vol 22 \endref







\ref \key 4
 \by  T. Kim
  \paper  On a $q$-analogue of the $p$-adic  log gamma functions
  and related integrals
 \jour  J. Number Theory
 \yr 1999
\pages 320--329 \vol 76 \endref




\ref \key 5
 \by  T. Kim
  \paper  $q$-Volkenborn integration
 \jour  Russ. J. Math. Phys.
 \yr 2002
\pages 288--299 \vol 9 \endref



\ref \key 6
 \by  T. Kim
  \paper An invariant  $p$-adic integral associated with Daehee numbers
 \jour Integral Transforms and Special Functions
 \yr 2002
\pages 65--69 \vol 13 \endref



\ref \key 7
 \by  T. Kim
  \paper Barnes-Euler multiple zeta functions
 \jour  Russ. J. Math. Phys.
 \yr 2003
\pages 185--192 \vol 6(2) \endref






\ref \key 8
 \by  T. Kim
  \paper  A note on multiple zeta functions
 \jour  JP Jour. Algebra Number Theory and Appl.
 \yr 2003
\pages 471--476 \vol 3(3) \endref






\ref \key 9 \by
 T\. Kim \paper
 Non-archimedean $q$-integrals associated with
multiple Changhee $q$-Bernoulli Polynomials \jour Russ. J. Math.
Phys. \vol 10 \yr 2003 \pages 91-98 \endref




\ref \key 10
 \by  T. Kim
  \paper  Remark on the multiple Bernoulli numbers
 \jour  Proc. Jangjeon Math. Soc.
 \yr 2003
\pages 185--192 \vol 6 \endref



\ref \key 11
 \by  T. Kim
  \paper   Sums of powers of consecutive $q$-integers
 \jour  Advan. Stud. Contem. Math.
 \yr 2004
\pages 15--18 \vol 9 \endref



\ref \key 12 \by
 T\. Kim \paper
 Analytic continuation of multiple $q$-zeta functions
and their values at negative integers \jour Russ. J. Math. Phys.
 \vol
11  \yr2004 \pages 71-76 \endref




\ref \key 13
 \by  T. Kim
  \paper  A note on multiple Dirichlet's $q$-$L$-function
 \jour  Advan. Stud. Contem. Math.
 \yr 2005
\pages 57--60 \vol 11 \endref





\ref \key 14 \by
 T. Kim  \paper
Power series and asymptotic series associated with the
$q$-analogue \ of two-variable $p$-adic $L$-function
 \jour Russ. J. Math. Phys.
 \yr 2005 \vol 12 \endref






\ref \key 15
 \by  T. Kim
  \paper  Multiple $p-$adic $L-$function
 \jour  Russ. J. Math. Phys.
 \yr 2006
\pages 151--157 \vol 13 \endref



\ref \key 16
 \by  T. Kim
  \paper   A new approach to $p$-adic $q$-$L$-functions
 \jour  Adv.
Stud. Contemp. Math. (Kyungshang)
 \yr 2006
\pages 61--72 \vol 12 \endref




\ref \key 16
 \by  T. Kim, L. C. Jang, S-H. Rim and H.-K. Pak
  \paper  On the twisted $q$-zeta functions and $q$-Bernoulli
  polynomials
 \jour  Far East J. Math.
 \yr 2003
\pages 13--21 \vol 13 \endref






\ref \key 18
 \by  N. Koblitz
  \paper A new proof of certain formulas for $p$-adic $L$-functions
 \jour  Duke Math. J.
 \yr 1979
\pages 455--468 \vol 40 \endref



\ref \key 19
 \by  J. Satho
  \paper $q$-analogue of Riemann's $\zeta$-function and $q$-Euler
  numbers
 \jour  J. Number Theory
 \yr 1989
\pages 346--362 \vol 31 \endref





\ref \key 20
 \by  M. Schlosser
  \paper $q$-analogue of the sums of  consecutive integers, squares, cubes, quarts, quints
   \jour Electro. J. Comb.
 \yr 2004
\pages  \vol 11 \# R71 \endref






\ref \key 21 \by K. Shiratani and S. Yamamoto
 \paper On a $p$-adic interpolation function for the Euler numbers and its
derivatives
 \jour Mem. Fac. Sci. Kyushu Univ. Ser. A \vol 39  \yr 1985 \pages 113-125  \endref




\ref \key 22
 \by  Y. Simsek; Yang, Sheldon
  \paper Transformation of four Titchmarch-type infinite integrals
  and generalized Dedekind sums associated with Lambert series
 \jour Adv. Stud. Contem. Math. (Kyungshang)
 \yr 2004
\pages  195--202 \vol  9 \endref



\ref \key 23
 \by  Y. Simsek
  \paper Theorems on twisted   $L$-functions and $q$-twisted Bernoulli numbers
 \jour Adv. Stud. Contem. Math.
 \yr 2005
\pages  205--218 \vol  11 \endref




\ref \key 24
 \by  Y. Simsek
  \paper Twisted $(h,q)-$Bernoulli numbers and polynomials related
  to $(h,q)-$zeta function and $L-$function
 \jour J. Math. Anal. Appl.(In press)
 \yr
\pages  \vol  \endref




\ref \key 25
 \by  Y. Simsek
  \paper  $q-$Dedekind type sums related to $q-$zeta function and
  basic $L-$series
 \jour J. Math. Anal. Appl.
 \yr 2006
\pages 333--351 \vol  318 \endref



\ref \key 26
 \by  Y. Simsek, D. KIm and S.-H. Rim
  \paper  On the two variable Dirichlet $q$-$L$-series
 \jour Adv. Stud. Contem. Math.
 \yr 2005
\pages 333--351 \vol  318 \endref




\ref \key 27
 \by  H. M.  Srivastava, T. Kim and Y. Simsek
  \paper  $q$-Bernoulli numbers and polynomials associated with
  multiple $q$-zeta functions and basic $L$-series
 \jour  Russ. J. Math. Phys.
 \yr 2005
\pages 241--268 \vol  12 \endref


\ref \key 28
 \by  Y. Simsek, D. Kim and S.-H. Rim
  \paper  A note on the sums of powes of consecutive $q$-integers
 \jour J. Appl. Funct. Different Equat.
 \yr 2006
\pages 63--70 \vol  1 \endref




\endRefs
\enddocument